\documentclass[a4paper,10pt]{article}
\usepackage{amsfonts}
\usepackage{amstext}
\usepackage{amsthm}
\usepackage{amsmath}
\usepackage{amssymb}
\usepackage{latexsym}
\usepackage[latin1]{inputenc}
\newcommand{\bl}{\hfill\rule{2mm}{2mm}}
\newcommand{\R}{\Bbb{R}}

\newtheorem{teor}{Theorem}

\newcommand{\n}{\noindent}

\begin{document}

\title{Equivalence of optimal $L^1$-inequalities on Riemannian Manifolds
\footnote{2010 Mathematics Subject Classification: 35A09, 35B44}
\footnote{Key words:    sharp Sobolev inequalities, best constant,
extremal maps}}

\author{\textbf{Jurandir Ceccon \footnote{\textit{E-mail addresses}:
ceccon@ufpr.br (J. Ceccon)}}\\ {\small\it Departamento de
Matem\'{a}tica, Universidade Federal do Paran\'{a},}\\
{\small\it Caixa Postal 019081, 81531-990, Curitiba, PR, Brazil}\\
\textbf{Leandro Cioletti \footnote{\textit{E-mail addresses}:
leandro.mat@gmail.com (L. Cioletti)}}\\ {\small\it Departamento de
Matem\'{a}tica, Universidade de Brasília,}\\
{\small\it UnB, 70910-900, Brasília, Brazil}} \maketitle

\markboth{abstract}{abstract}
\addcontentsline{toc}{chapter}{abstract}

\hrule \vspace{0,2cm}

\n {\bf Abstract}

Let $(M,g)$ be a smooth compact Riemannian manifold of dimension
$n \geq 2$. This paper concerns to the validity of the optimal
Riemannian $L^1$-Entropy inequality
\[
{\bf Ent}_{dv_g}(u) \leq n \log \left( A_{opt} \|D u\|_{BV(M)} +
B_{opt}\right)
\]
\n for all $u \in BV(M)$ with $\|u\|_{L^1(M)} = 1$ and existence
of extremal functions. In particular, we prove that this optimal
inequality is equivalent a optimal $L^1$-Sobolev inequality
obtained by Druet \cite{Druet}.

\vspace{0,5cm} \hrule\vspace{0.2cm}

\section{Introduction}

The isoperimetric problem on the Euclidean space $\mathbb{R}^n$
consist in finding among all the domains with a given fixed 
volume one that has the lowest surface area.
The solution in this case is a sphere.
This property is precisely expressed in terms of the 
Isoperimetric inequality for domains in $\mathbb{R}^n$,
that is, if $\Omega$ is a domain with volume $|\Omega|$
and surface area $|\partial \Omega|$ then
\begin{equation}\label{di}
\frac{|\partial \Omega|}{|\Omega|^{\frac{n - 1}{n}}} \geq
\frac{1}{K(n,1)} \; ,
\end{equation}
where
$
K(n,1)^{-1} = n^{(1-1/n)} (\omega_{n-1})^{\frac{1}{n}}
$
and $\omega_{n-1}$ denotes the volume of the unit ball in the Euclidean
space $\mathbb{R}^{n-1}$.
The equality is attained iff $\Omega$ is a sphere.
Observe that the last statement implies that $K(n,1)$ is the best constant 
for the inequality \eqref{di}. 
The Isoperimetric inequality shows up on other branchs of mathematics. 
For example, it has been used to prove the non-uniqueness of the Gibbs measures 
of the Ising model on the lattice $\mathbb{Z}^n\equiv \mathbb{Z}\times\ldots\times\mathbb{Z}$.
In this context, the inequality is generalized as follows. 
We consider $\mathbb{Z}^n$ as a metric space, where the distance between $x,y\in\mathbb{Z}^n$
is defined in terms of the $\ell_1$ norm. 
We also look at $\mathbb{Z}^n$ as a graph $(\mathbb{Z}^n,\mathbb{E}^n)$, 
where $\mathbb{Z}^n$ is the vertex set and 
$\mathbb{E}^n\equiv \{\{x,y\}\in \mathbb{Z}^n\times\mathbb{Z}^n; |x-y|_1=1\}$ is the edge set.
The discrete Isoperimetric inequality says that for any
fixed integer $n\geq 2$ and any finite subset $\Omega\subset\mathbb{Z}^n$, we have that 
\[
\frac{|\partial \Omega|}{|\Omega|^{\frac{n-1}{n}}}
	\geq \frac{1}{2n},
\]
where 
$	\partial \Omega
	=
	\{\{i,j\}\in\ \mathbb{E}^d: i\in\Omega, j\in\Omega^c\}
$
and $|\Omega|$ and $|\partial \Omega|$
denotes the cardinality of $\Omega$ and $\partial \Omega$, respectively. 
Note that $2n$ is exactly the volume of the unit 
sphere on the $\ell_1$ norm. 
It is worth pointing out that the proof of this discrete inequality
is similar on spirit of our proof on the continuous setting 
and is based on the entropy inequalities. 
It is possible that the equivalences obtained here can 
be extended to the discrete case but we will not develop this point here. More information about this connection can be found in \cite{sallof-coste2}. For more details about the Isoperimetric inequality see the excellent Osserman's work
\cite{Ro} and references therein.

Now we consider a more analytic context.
We say that a
function $u \in L^1(\R^n)$ has {\bf bounded variation} if

\[
\|D u\|_{BV(\R^n)} = \sup \left\{\int_{\R^n} u\,
\mbox{div}(\varphi)\, dx;\, \varphi \in C_0^1(\R^n,\R^n),
|\varphi| \leq 1 \right\} < \infty \; .
\]
\n The space of all bounded variation functions is denoted by
$BV(\R^n)$. The optimal Euclidean $L^1$-Sobolev inequality in
$BV(\R^n)$ states that for all $u \in BV(\R^n)$ we have
\begin{equation}\label{dos}
\|u\|_{L^{1^*}(\R^n)} \leq K_0 \|D u\|_{BV(\R^n)} \; ,
\end{equation}
where $1^* = \frac{n}{n - 1}$ is the critical Sobolev exponent and

\[
K_0^{-1} = \inf\{ \|Du\|_{BV(\R^n)}; u \in BV(\R^n),
\|u\|_{L^{1^*}(\R^n)} = 1 \} \;
\]
\n is the best constant for this inequality. This inequality it
was studied by Federer and Fleming \cite{ff}, Fleming and Rishel
\cite{fr} and Maz$'$ja \cite{VM}. In this case, the characteristic
functions of the balls are extremal functions for the optimal
$L^1$-Sobolev inequality and explicit value for best constant is
given by
\[
K_0 = K(n,1) \; .
\]
This inequalities gains in interest if we realize that the geometric inequality \eqref{di}
and the analytic inequality \eqref{dos} are equivalents.
This relation was pointed
out independently by Federer and Fleming \cite{ff} and Maz$'$ja
\cite{VM}. Let's move on to the manifold setting.

Let $(M,g)$ be a smooth compact Riemannian manifold without
boundary of dimension $n \geq 2$. We say that $u \in L^1(M)$ is a
bounded variation function if
\[
\|D_g u\|_{BV} =\sup
    \left\{
        \int_M u\, \mbox{div}_g(X)\, dv_g; X \in
        \Gamma(TM), |X|_g(x) \leq 1 \hspace{0,2cm} \mbox{for all}
        \hspace{0,2cm} x \in M
    \right\}
< \infty \; ,
\]

\n where $\Gamma(TM)$ is the set of all vector fields over $M$
with divergent in $n$-fold Cartesian product $C^1(M) \times \cdots
\times C^1(M)$. We denote by $BV(M)$ the space of bounded
variation functions. The optimal Riemannian $L^1$-Sobolev
inequality was obtained by Druet \cite{Druet}. He proved that for
all $u \in BV(M)$ the following inequality holds
\begin{equation}\label{mdos}
\|u\|_{L^{1^*}(M)} \leq K(n,1) \|D_g u\|_{BV(M)} + B(1)
\|u\|_{L^1(M)} \; ,
\end{equation}

\n where $K(n,1)$ is the {\bf first best constant} and $B(1)$ is
the {\bf second best constant} for the optimal Riemannian
$L^1$-Sobolev inequality.

Moreover, Druet also proved that the only extremal 
functions for \eqref{mdos} are the characteristic 
function of some $\Omega_0\in \Sigma$, where 
\(
\Sigma = \{\Omega \subset M; \chi_{\Omega} \in BV(M)\} \; .
\)
Finally, Druet observed that for $\Omega \in \Sigma$, 
the optimal inequality \eqref{mdos} is equivalent to the 
following isoperimetric inequality
\begin{equation}\label{isor}
|\Omega|^{\frac{n - 1}{n}} \leq K(n,1) |\partial \Omega| + B(1)
|\Omega|,
\end{equation}

\n where $|\Omega|$ and $|\partial \Omega|$ are the Riemannian
volume and area of $\Omega$ and $\partial \Omega$, respectively
and $B(1)$ is the best constant in this inequality. 
For more details about this equivalence and additional
references see Druet \cite{Druet}.

The aim of this work is to show that these inequalities 
are also equivalent to the Entropy and Gagliardo-Nirenberg
inequalities on both context Euclidean and Riemannian, and determine 
the equivalence among its extremal functions.

\section{Euclidean inequality}

For any function $u \in BV(\R^n)$ we have that its entropy, 
with respect to the Lebesgue measure, is well defined and given
by the following expression
\[
{\bf Ent}_{dx}(u) = \int_{\R^n} |u| \log|u| dx \; .
\]
\n The optimal Euclidean {\bf entropy inequality} states that
for all $u \in BV(\R^n)$ with $\|u\|_{L^1(\R^n)} = 1$,
\begin{equation}\label{dle}
{\bf Ent}_{dx}(u) \leq n \log \left( L(n,1) \|\nabla
u\|_{BV(\R^n)} \right) \; ,
\end{equation}
where $L(n,1)$ is the best constant for this inequality. 
Note that the existence of this best constant is guaranteed  
by the same argument we use in \eqref{dos}.
This inequality it was studied by Beckner \cite{Be} and Ledoux
\cite{Le}. For the optimal entropy inequality, Beckner showed that
the extremal functions for (\ref{dle}) are the characteristic
functions of the Euclidean balls $B(x_0,r) = \{x\in \R^n; \|x -
x_0\| < r\}$, with $|B(x_0,r)| = 1$. So in this case, it is
easy to check that
\begin{equation}\label{sl}
K(n,1) = L(n,1) \; .
\end{equation}
In addition, if $\Omega$ is a domain in $\mathbb{R}^n$
and $\lambda>0$ is chosen so that
\[
\lambda \int_M \chi_{\Omega} dx = 1
\]
then by using $\lambda \chi_{\Omega}$ in \eqref{dle}, we
obtain the Isoperimetric inequality \eqref{di}.

Now we show that the entropy inequality can be obtained 
by the Gagliardo-Nirenberg inequality. Let $1 \leq q <
r < 1^*$. By using the interpolation inequality and the optimal
Euclidean $L^1$-Sobolev inequality, we obtain the Euclidean
Gagliardo-Nirenberg inequality, which says that for any function
$u \in BV(\R^n)$ we have
\begin{equation}\label{dgne}
\|u\|_{L^r(\R^n)}^{\frac{1}{\theta}} \leq K(n,1) \|D
u\|_{BV(\R^n)} \|u\|_{L^q(\R^n)}^{\frac{1 - \theta}{\theta}}\, ,
\end{equation}

\n where $\theta = \frac{n(r - q)}{r(q(1 - n) + n)} \in (0,1)$.
Let
\[
A(n,q,r)^{-1} = \inf \left\{
            \|Du\|_{BV(\R^n)} \|u\|_{L^q(\R^n)} ^{\frac{1 -
            \theta}{\theta}}; u\in BV(\R^n) , \|u\|_{L^r(\R^n)} = 1
\right\}
\]
be the best possible constant for this inequality.
If $\chi_{B(0,r)}$ denotes the characteristic function of the Euclidean ball
of radius $r>0$, then we have that $\chi_{B(0,r)} \in BV(\R^n)$ and an
easy computation shows that
\[
\|\chi_{B(0,r)}\|_{L^r(\R^n)}^{\frac{1}{\theta}}
= K(n,1)
\|D \chi_{B(0,r)}\|_{BV(\R^n)}
\cdot\|\chi_{B(0,r)}\|_{L^q(\R^n)}^{\frac{1-\theta}{\theta}} \; .
\]
Note that the above equality actually proved that
\begin{equation}\label{equa}
A(n,q,r) = K(n,1) \; ,
\end{equation}
for all $1 \leq q < r < 1^*$.  Therefore the inequality
(\ref{dgne}) is in fact the {\bf optimal Euclidean
$L^1$-Gagliardo-Nirenberg inequality} and characteristic function
of the balls are extremal functions for \eqref{dgne}.

Proceeding, with minor modifications, as in \cite{Ba} 
(see also Section 3) we can verify that the optimal
Euclidean $L^1$-Gagliardo-Nirenberg inequality implies 
the optimal Euclidean $L^1$-entropy inequality \eqref{dle}.
Piecing together these information we conclude that 
\eqref{di} $\Rightarrow$ \eqref{dos} $\Rightarrow$ \eqref{dgne} $\Rightarrow$
\eqref{dle} $\Rightarrow$ \eqref{di}.

Before proceed, we remark that the equality \eqref{equa} is the
key point in Section 3 to prove the main result of this paper.
\section{The Riemannian inequality and the main result}

Using the interpolation inequality and
the optimal Riemannian $L^1$-Sobolev inequality we get
for any function $u \in BV(M)$, $1 \leq q < r < 1^*$ that
\begin{equation}\label{dognr}
\|u\|_{L^r(M)}^{\frac{1}{\theta}} \leq \left( K(n,1) \|D_g
u\|_{BV(M)} + B(1) \|u\|_{L^1(M)}\right) \|u\|_{L^q(M)}^{\frac{1 -
\theta}{\theta}}\, , \tag{$I_{q,r}(K(n,1),B(1))$}
\end{equation}
where $\theta = \frac{n(r - q)}{r(q(1 - n) + n)}$ is the
interpolation parameter.
The {\bf first Riemannian $L^1$-Gagliardo-Nirenberg best constant}
is defined by
\[
A_{opt} = \inf \{ A \in \R:\; \mbox{there exists} \hspace{0,18cm}
B \in \R \hspace{0,18cm} \mbox{such that} \hspace{0,18 cm}
I_{q,r}(A,B) \hspace{0,18cm} \mbox{is valid}\}\, .
\]

\n
Using a partition unity and a similar argument as in \cite{DHV}
together with (\ref{equa}), we can verify that
\[
{\cal A}_{opt} = A(n,q,r) = K(n,1) \; ,
\]
\n for all $1 \leq q < r < 1^*$. So this equality shows that the
{first \bf optimal Riemannian $L^1$-Gagliardo-Nirenberg
inequality} (\ref{dognr}) is valid for all $u \in BV(M)$ and
$K(n,1)$ is the {\bf first best constant}.

It follows from \cite{Druet} that every extremal function for the
Sobolev inequality is of the form
\(
u_0 = \lambda \chi_{\Omega_0} \; ,
\)
\n for some $\lambda \in \R$ and $\Omega_0 \in \Sigma$.
We see at once that for such functions the equality
in $I_{q,r}(K(n,1),B(1))$ is verified.
Consequently the {\bf second Riemannian $L^1$-Gagliardo-Nirenberg best
constant} is given by
\[
B(1) = \inf \{ B \in \R; I_{q,r}(K(n,1),B) \hspace{0,2cm} \mbox{is
valid}\} \; .
\]
So we have that {\bf optimal Riemannian
$L^1$-Gagliardo-Nirenberg inequality}
\begin{equation}\label{dro}
\|u\|_{L^r(M)}^{\frac{1}{\theta}} \leq \left( K(n,1) \|D_g
u\|_{BV(M)} + B(1) \|u\|_{L^1(M)}\right) \|u\|_{L^q(M)}^{\frac{1 -
\theta}{\theta}}\, ,
\end{equation}
is valid for all $u \in BV(M)$ and $1 \leq q < r < 1^*$.
Notice that the Sobolev extremal functions are also
Gagliardo-Nirenberg extremal functions.

In the sequel we show how to use the inequality (\ref{dro}) (with $q
= 1$) to obtain the optimal Riemannian $L^1$-entropy inequality.
The proof is based on the Bakry, Coulhon, Ledoux and Sallof-Coste
argument given in \cite{Ba}. Consider the optimal Riemannian
$L^1$-Gagliardo-Nirenberg inequality

\[
\|u\|_{L^r(M)}^{\frac{1}{\theta}} \left(\leq K(n,1) \|D_g
u\|_{BV(M)} + B(1) \|u\|_{L^1(M)} \right) \|u\|_{L^1(M)}^{\frac{1
- \theta}{\theta}} \; ,
\]

\n for all $u \in BV(M)$. By taking the logarithm on both sides above
and use the definition of $\theta$, we get
\[
\frac{r(1 - n + n)}{n} \frac{1}{r - 1} \log \left(
\frac{\|u\|_{L^r(M)}}{\|u\|_{L^1(M)}} \right) \leq \log \left(
K(n,1) \frac{\|D_g u\|_{BV(M)}} {\|u\|_{L^1(M)}} + B(1)\right).
\]
Taking the limit when $r\to 1^+$, on the above expression, we obtain
\[
\frac{1}{n} \lim_{r \rightarrow 1^+} \frac{1}{r - 1} \log
\left(\frac{\|u\|_{L^r(\R^n)}}{\|u\|_{L^1(\R^n)}} \right)\leq \log
\left(K(n,1) \frac{\|D_g u\|_{L(M)}}{\|u\|_{L^1(M)}} + B(1)\right)
\; .
\]
\n
To evaluate the remainder limit, we first observe that
\begin{align*}
\log \left(\frac{\|u\|_{L^r(M)}}{\|u\|_{L^1(M)}} \right) &=
\frac{1}{r} \log(\|u\|_{L^r(M)}^r) -  \log(\|u\|_{L^1(M)})
\\[0.2cm]
&= \frac{1 - r}{r} \log(\|u\|_{L^1(M)}) + \frac{1}{r} \left(
\log(\|u\|_{L^r(M)}^r) - \log(\|u\|_{L^1(M)}) \right).
\end{align*}
\n
Next, we apply two times the mean value theorem, obtaining
\[
\lim_{r \rightarrow 1^+} \frac{1}{r - 1} \log
\left(\frac{\|u\|_{L^r(M)}}{\|u\|_{L^1(M)}} \right) = \int_{M}
\frac{|u|}{\|u\|_{L^1(M)}} \log
\left(\frac{|u|}{\|u\|_{L^1(M)}}\right) dx.
\]

\n From the above equation it follows that
\[
\int_{M} |u| \log (|u|) dv_g \leq n \log \left(K(n,1) \|D_g
u\|_{BV(M)} + B(1) \right)\ ,
\]
\n for all $u \in BV(M)$ with $\|u\|_{L^1(M)} = 1$.
As in the previous section we define
\[
{\bf Ent}_{dv_g}(u) = \int_M |u| \log |u| dv_g \; .
\]
\n As a consequence of the previous inequality we have that
\begin{equation}\label{ent}
{\bf Ent}_{dv_g}(u) \leq n \log \left(L(n,1) \|D_g u\|_{BV(M)} +
B(1) \right) \; , \tag{$Ent(L(n,1),B(1)$}
\end{equation}
\n for all $u \in BV(M)$ with $\|u\|_{L^1(M)} = 1$. We shall remember
that $K(n,1) = L(n,1)$.
\n Now we consider the {\bf optimal Riemannian
$L^1$-entropy inequality}
\[
{\bf Ent}_{dv_g}(u) \leq n \log \left(L_{opt} \|D_g u\|_{BV(M)} +
B \right)
\]
where $u \in BV(M)$ with $\|u\|_{L^1(M)} = 1$,
$B\in\mathbb{R}$ and the {\bf first Riemannian $L^1$-entropy
best constant} is defined by
\[
L_{opt} = \inf \{ A \in \R:\; \mbox{there exists} \hspace{0,18cm}
B \in \R \hspace{0,18cm} \mbox{such that} \hspace{0,18 cm}
Ent(A,B) \hspace{0,18cm} \mbox{is valid}\}\, .
\]
From the definition of best constant and the validity of
$Ent(L(n,1),B(1))$ follows that $L_{opt} \leq L(n,1)$.
On the other hand, the equality between these two constants
requires a proof.

Assuming that $Ent(L_{opt},B)$ is valid for some $B \in \R$,
one can define the {\bf second Riemannian $L^1$-Entropy best
constant} by
\[
B_{opt} = \inf \{ B \in \R; Ent(L_{opt},B) \hspace{0,2cm} \mbox{is
valid}\} \; .
\]
Our main result states that $Ent(L_{opt},B_{opt})$ is valid,
moreover the Riemannian first best constant is equals to the Euclidean first best constant.
We also prove that the second best constant for Sobolev and entropy inequalities
are the same.

\begin{teor} \label{teore}
Let $(M,g)$ be a smooth compact Riemannian manifold without
boundary of dimension $n \geq 2$. Then $Ent(L_{opt},B(1))$ is
valid. In addition, $L_{opt} = L(n,1)$ and $B_{opt} = B(1)$.
\end{teor}

\noindent {\bf Remark.} Together with the results in \cite{Druet}
the Theorem \ref{teore} can be used to complete the equivalence
of the four $L^1$ optimal inequalities considered on this work.

Indeed, in \cite{Druet} we have that isoperimetric $\Rightarrow$
Sobolev. In the Section 3 we shown that Sobolev $\Rightarrow$
Gagliardo-Nirenberg and Gagliardo-Nirenberg + Theorem 1
$\Rightarrow$ entropy. Finally, by taking $\Omega \in \Sigma$
and $\lambda > 0$ such that
\[
\int_M \lambda \; \chi_{\Omega} dv_g = 1 \; ,
\]
and replacing $\lambda \; \chi_{\Omega}$ in $Ent(L(n,1),B(1))$
we get the optimal isoperimetric inequality. 

So the four $L^1$-optimal inequalities considered here 
are all equivalent and we have that \eqref{isor} $\Rightarrow$ \eqref{mdos}
$\Rightarrow$ \eqref{dognr} $\Rightarrow$ \eqref{ent}
$\Rightarrow$ \eqref{isor} .

\n {\bf Proof of Theorem \ref{teore}:} We proceed to show that
\[
L_{opt} = L(n,1) \;.
\]
From this equality it may be concluded that $Ent(L(n,1),B(1))$ is the best
$L^1$-entropy.
To prove the above equality, we use the definition of
$L_{opt}$. That is, by definition, for all $s > 0$ there exists
$B(s)$ such that
\[
\int_M |u| \log (|u|) dv_g \leq n \log \left((L_{opt} + s) \|D_g
u\|_{BV(M)} + B(s) \right)\ ,
\]
\n for all $u \in BV(M)$ with $\|u\|_{L^1(M)} = 1$. This is equivalent
to the following inequality

\begin{equation} \label{EF}
\frac{1}{\|u\|_{L^1(M)}} \int_M |u| \log(|u|) dv_g + (n - 1) \log
(\|u\|_{L^1(M)}) \leq n \log \left(  (L_{opt} + s) \|D_g
u\|_{BV(M)} + B(s) \int_M |u| dv_g \right)
\end{equation}

\n for all $u \in BV(M)$. Let $r > 0$ be such that

\[
\int_{B(0,r)}  dx = 1\; ,
\]

\n where $B(0,r) = \{x\in\R^n; \|x\| < r\}$ is the Euclidean ball.
For $\varepsilon > 0$, we define
$\chi_{\varepsilon} = (r\varepsilon)^{-n} \chi_{B_{g}(0,\varepsilon)}$, where
$\chi_{B_g(0,\varepsilon)}$ is the characteristic function of the
Riemannian ball $B_g(0,\varepsilon) = \{x \in M; d_g(x_0,x) \leq
\varepsilon \}$ and $x_0 \in M$.
In the normal coordinates in around $x_0$, we have
\[
\int_M \chi_\varepsilon dv_g = (\varepsilon r )^{-n}
\int_{B(x_0,\varepsilon)} dv_g = \int_{B(0,r)} dv_{g_\varepsilon}
\]
and $dv_{g_\varepsilon} = \sqrt{(g_{ij}(\exp_{x_0}(\varepsilon
x)))}dx \rightarrow dx$ when $\varepsilon \rightarrow 0$, where
$g_{ij}$ are the coefficients of the metric $g$ in the normal
coordinates. We also have

\[
\|D_g \chi_{\varepsilon}\|_{BV(M)} = (\varepsilon r)^{-n}
\int_{\partial B(x_0,\varepsilon)} dv_g = (\varepsilon r)^{-1}
\int_{\partial B(0,r)} dv_{g_\varepsilon} \; .
\]

\n By replacing this identity in (\ref{EF}) we get
\[
-n \log (\varepsilon r) + (n - 1) \log |B(0,r)|_{g_\varepsilon}
\leq -n \log (\varepsilon r) + n \log\left( (L_{opt} + s)
|\partial B(0,r)|_{g_\varepsilon} + \varepsilon r B(s)
|B(0,r)|_{g_\varepsilon}\right) \; .
\]

\n Because $|B(0,r)|_{\xi} = 1$ ($\xi$ denotes the Euclidean
metric) and by taking the limits when $\varepsilon \rightarrow 0$ and
$s \rightarrow 0$ in this order, we obtain
\[
0 \leq n \log(L_{opt} |\partial B(0,r)|_{\xi}) \; .
\]
By the choice of $r$ and remembering again that the characteristic
function of the Euclidean ball $B(0,r)$ is the extremal function
in (\ref{dle}), we have now reached the following identity
\[
0 = n \log( L(n,1) |\partial B(0,r)|_{\xi}) \; .
\]
\n
From where we conclude that $L(n,1) \leq L_{opt}$, which
proves the desired equality. Thus (\ref{ent}) is the
optimal Riemannian $L^1$-entropy.

Now we compute the second best constant.
Consider $\Omega_0 \in
\Sigma$ and $k > 0$, such that
\[
\int_M k \chi_{\Omega_0} dv_g = 1 \;,
\]
such that, $\chi_{\Omega_0}$ is the extremal function in the
optimal Sobolev inequality (\ref{dos}).
We can immediately check that
the optimal inequality $Ent(L(n,1),B(1))$ becomes an equality when
we evaluate in the function $k \chi_{\Omega_0}$.
Therefore $B(1) = B_{opt}$.
\bl

\section{Equivalence between extremal functions}

As we already observed, the extremal functions for the optimal
$L^1$-Sobolev inequality (\ref{mdos}) are characteristic
functions. We also remarked that these functions are also extremal
for the optimal $L^1$-Gagliardo-Nirenberg inequality
(\ref{dognr}). We can also prove, using the limit process employed
in the Section 3, that the extremal functions for the optimal
$L^1$-Gagliardo-Nirenberg \eqref{dognr} are the extremal functions
for the optimal $L^1$-Entropy inequality \eqref{ent}.

If we prove that the extremal functions for the optimal
$L^1$-Entropy inequality (\ref{ent}) are the extremal functions
for the optimal $L^1$-Sobolev inequality (\ref{mdos}), we have
that the set of the extremal functions for the four optimal $L^1$
inequalities considered here are the same.

\n
{\bf Claim.} If $u_0$ is an extremal function for (\ref{ent})
$\Longrightarrow$ $u_0$ is an extremal function for (\ref{dos}).
In fact, from the Jensen Inequality for any $u_0$ such that $\|u_0\|_{L^1(M)} = 1$
we get
\[
\log \int_M |u_0|^{1^*} dv_g = \log \int_M |u_0|^{1^* - 1} |u_0|
dv_g \geq \int_M \log(|u_0|^{\frac{1}{n - 1}}) |u_0| dv_g =
\frac{1}{n - 1} \int_M \log(|u_0|) |u_0| dv_g.
\]
By using (\ref{ent}) it follows that
\[
\log \int_M |u|^{1^*} dv_g \geq 1^* \log \left(K(n,1)
\|Du_0\|_{BV(M)} + B(1) \right),
\]
that is,
\[
\|u_0\|_{L^{1^*}(M)} \geq K(n,1) \|Du_0\|_{BV(M)} + B(1).
\]
This inequality shows that $u_0$ is an extremal function for
(\ref{mdos}). \bl

{\bf Acknowledgments.} The first author was partially supported
by CAPES through INCTmat and second author is supported by
FEMAT.

\end{document}